# On the Differential Operators with Periodic Matrix Coefficients

O. A. Veliev

Depart. of Math., Dogus University, Acıbadem, Kadıköy,

Istanbul, Turkey. e-mail: oveliev@dogus.edu.tr

**Abstract**

In this article we obtain asymptotic formulas for eigenvalues and eigenfunctions of the operator generated by a system of ordinary differential equations with summable coefficients and quasiperiodic boundary conditions. Then using these asymptotic formulas, we find conditions on the coefficients for which the number of gaps in the spectrum of the self-adjoint differential operator with the periodic matrix coefficients is finite.

Key Words: Differential operator, Spectrum, Asymptotic formulas.

AMS Mathematics Subject Classification: 34L05, 34L20.

Let $L(P_2, P_3..., P_n)$ be the differential operator generated in the space $L_2^m(-\infty, \infty)$ of vector valued functions by the differential expression

$$(-i)^n y^{(n)}(x) + (-i)^{n-2} P_2(x) y^{(n-2)}(x) + \sum_{\nu=3}^{n} P_\nu(x) y^{(n-\nu)}(x), \qquad (1)$$

where $n \geq 2$, $P_\nu(x) = (p_{\nu,i,j}(x))$ is a $m \times m$ matrix with the complex-valued summable entries $p_{\nu,i,j}(x)$, $P_\nu(x+1) = P_\nu(x)$ for $\nu = 2, 3, ...n$. It is well-known that ( see [1, 2, 4]) the spectrum of the operator $L(P_2, P_3..., P_n)$ is the union of the spectra of the operators $L_t(P_2, P_3..., P_n)$ for $t \in [0, 2\pi)$ generated in $L_2^m(0,1)$ by the expression (1) and the quasiperiodic conditions

$$U_\nu(y) \equiv y^{(\nu)}(1) - e^{it} y^{(\nu)}(0) = 0, \ \nu = 0, 1, ..., (n-1). \qquad (2)$$

Note that $L_2^m(a,b)$ is the set of the vector valued functions $f = (f_1, f_2, ..., f_m)$ with $f_k \in L_2(a,b)$ for $k = 1, 2, ..., m$. The norm $\|.\|$ and inner product $(.,.)$ in $L_2^m(a,b)$ are defined by

$$\|f\|^2 = \left( \int_a^b |f(x)|^2 dx \right), \ (f,g) = \int_a^b \langle f(x), g(x) \rangle dx,$$

where $|.|$ and $\langle ., . \rangle$ are the norm and inner product in $\mathbb{C}^m$.

Let us introduce some preliminary results and describe the results of this paper. Clearly,

$$\varphi_{k,1,t} = \left( e^{i(2\pi k + t)x}, 0, ..., 0 \right), \ \varphi_{k,2,t} = \left( 0, e^{i(2\pi k + t)x}, 0, ..., 0 \right), ..., \varphi_{k,m,t} = \left( 0, 0, ..., 0, e^{i(2\pi k + t)x} \right)$$

are the eigenfunctions of the operator $L_t(0)$ corresponding to the eigenvalue $(2\pi k + t)^n$, where $k \in \mathbb{Z}$, and the operator $L_t(P_2, ..., P_n)$ is denoted by $L_t(0)$ when

$P_2(x) = 0, ..., P_n(x) = 0$. Furthermore, for brevity of notation, the operators $L_t(P_2, ..., P_n)$ and $L(P_2, ..., P_n)$ are denoted by $L_t(P)$ and $L(P)$ respectively. It easily follows from the classical investigations [5, chapter 3, theorem 2] that the large eigenvalues of the operator $L_t(P)$ consist of $m$ sequences

$$\{\lambda_{k,1}(t) : \mid k \mid \geq N\},\ \{\lambda_{k,2}(t) :\mid k \mid \geq N\}, ...,\ \{\lambda_{k,m}(t) :\mid k \mid \geq N\} \tag{3}$$

satisfying the following, uniform with respect to $t$ in $[0, 2\pi)$, asymptotic formulas

$$\lambda_{k,j}(t) = (2\pi k + t)^n + O\left(k^{n-1-\frac{1}{2m}}\right) \tag{4}$$

for $j = 1, 2, ..., m$, where $N$ is a sufficiently large positive number, that is $N \gg 1$. We say that the formula $f(k,t) = O(h(k))$ is uniform with respect to $t$ in a set $S$ if there exists a positive constant $c_1$, independent of $t$, such that $\mid f(k,t)) \mid < c_1 \mid h(k) \mid$ for all $t \in S$ and $k \in \mathbb{Z}$. Thus the formula (4) means that there exists a positive numbers $N$ and $c_1$, independent of $t$, such that

$$\mid \lambda_{k,j}(t) - (2\pi k + t)^n \mid < c_1 \mid k \mid^{n-1-\frac{1}{2m}},\ \forall \mid k \mid \geq N,\ \forall t \in [0, 2\pi). \tag{5}$$

Note that in the classical investigations (see [5]) in order to obtain the asymptotic formulas of high accuracy, by using the classical asymptotic expansions for solutions of the matrix equation

$$(-i)^n Y^{(n)} + (-i)^{n-2} P_2 Y^{(n-2)} + \sum_{\nu=3}^{n} P_\nu Y^{(n-\nu)} = \lambda Y, \tag{6}$$

it is required that the coefficients be differentiable. Thus these classical methods never permits us to obtain the asymptotic formulas of high accuracy for the operator $L_t(P)$ with nondifferentiable coefficients. However, the method suggested in this paper is independent of smoothness of the coefficients. In this paper, by the suggested method, we obtain the uniform asymptotic formulas of high accuracy for the eigenvalues $\lambda_{k,j}(t)$ and for the corresponding normalized eigenfunctions $\Psi_{k,j,t}(x)$ of $L_t(P)$ when the entries $p_{2,i,j}(x)$, $p_{3,i,j}(x), ..., p_{n,i,j}(x)$ of $P_2(x), P_3(x), ..., P_n(x)$ belong to $L_1[0,1]$, that is, when there is not any condition about smoothness of the coefficients. Then using these formulas, we find the conditions on the coefficients for which the number of the gaps in the spectrum of the self-adjoint differential operator $L(P)$ is finite.

Now let us describe the scheme of the paper. The inequality (5) shows that the eigenvalue $\lambda_{k,j}(t)$ of $L_t(P)$ is close to the eigenvalue $(2k\pi + t)^n$ of $L_t(0)$. To analyze the distance of the eigenvalue $\lambda_{k,j}(t)$ of $L_t(P)$ from the other eigenvalues $(2p\pi + t)^n$ of $L_t(0)$, which is important in perturbation theory, we take into account the following situations. If the order $n$ of the differential expression (1) is odd number, $n = 2r - 1$, and $\mid k \mid \gg 1$, then the eigenvalue $(2\pi k + t)^n$ of $L_t(0)$ lies far from the other eigenvalues $(2p\pi + t)^n$ of $L_t(0)$ for all values of $t \in [0, 2\pi)$. We have the same situation if $n = 2r$ and $t$ does not lie in the small neighborhoods of $0$ and $\pi$. However, if $n$ is even number and $t$ lies in the neighborhoods of $0$ and $\pi$, then the eigenvalue $(2\pi k + t)^n$ is close to the eigenvalues $(2\pi(-k) + t)^n$ and $(2\pi(-k-1) + t)^n$ respectively. For this reason instead of $[0, 2\pi)$ we consider $t \in [-\frac{\pi}{2}, \frac{3\pi}{2})$ and use the following notation

**Notation 1** *Case 1 :* (a) $n = 2r - 1$ *and* $t \in [-\frac{\pi}{2}, \frac{3\pi}{2})$, (b) $n = 2r$ *and* $t \in T(k)$, *where*
$T(k) = [-\frac{\pi}{2}, \frac{3\pi}{2}) \backslash ((-\frac{1}{\ln|k|}, \frac{1}{\ln|k|}) \cup (\pi - \frac{1}{\ln|k|}, \pi + \frac{1}{\ln|k|})).$

*Case 2:* $n = 2r$ *and* $t \in ((-(\ln|k|)^{-1}, (\ln|k|)^{-1}))$.

*Case 3:* $n = 2r$ *and* $t \in (\pi - (\ln|k|)^{-1}, \pi + (\ln|k|)^{-1}))$.

*Denote by* $A(k, n, t)$ *the sets* $\{k\}, \{k, -k\}, \{k, -k-1\}$ *for the Cases 1, 2, 3 respectively.*

By (5) there exists a positive constant $c_2$, independent of $t$, such that the inequalities

$$|(2k\pi + t)^n - (2\pi p + t)^n| > c_2(\ln|k|)^{-1}(||k| - |p|| + 1)(|k| + |p|)^{n-1},$$
$$|\lambda_{k,j}(t) - (2\pi p + t)^2| > c_2(\ln|k|)^{-1}(||k| - |p|| + 1)(|k| + |p|)^{n-1}, \quad (7)$$

where $|k| > N$, hold in Case 1, Case 2, and Case 3 for $p \neq k$, for $p \neq k, -k$, and for $p \neq k, -(k+1)$ respectively. To avoid the listing of these cases, using the Notation 1, we see that the inequalities in (7) hold for $p \notin A(k,n,t)$. Using this and taking into account that if $p \neq \pm k, \pm(k+1), \pm(k+2)$, then the inequalities obtained from (7) by deleting the multiplicand $(\ln|k|)^{-1}$ in the right-hand side of (7) hold, we obtain the following useful formulas

$$\sum_{p:p>d} \frac{p^{n-2}}{|\lambda_{k,j}(t) - (2\pi p + t)^n|} = O(\frac{1}{d}), \quad (8)$$

$$\sum_{p:p \notin A(k,n,t)} \frac{p^{n-2}}{|\lambda_{k,j}(t) - (2\pi p + t)^n|} = O(\frac{\ln|k|}{k}), \quad (9)$$

$$\sum_{p:p \notin A(k,n,t)} \frac{k^{2n-4}}{|\lambda_{k,j}(t) - (2\pi p + t)^n|^2} = O(\frac{(\ln|k|)^2}{k^2}), \quad (10)$$

$$\sum_{p:p \notin A(k,n,t)} \frac{p^{2n-4}}{|\lambda_{p,j}(t) - (2\pi k + t)^n|^2} = O(\frac{(\ln|k|)^2}{k^2}), \quad (11)$$

where $d > 2|k|$, $|k| \geq N$, $N$ is defined in (5), and the formulas (8)-(11) are uniform with respect to $t$ in $[-\frac{\pi}{2}, \frac{3\pi}{2})$.

To obtain the asymptotic formulas we use (8)-(11) and consider the operator $L_t(P)$ as perturbation of $L_t(C)$ by $L_t(P) - L_t(C)$, where $C = \int_0^1 P_2(x)\,dx$, $L_t(C)$ is the operator generated by (2) and by the expression

$$(-i)^n y^{(n)}(x) + (-i)^{n-2} C y^{(n-2)}(x). \quad (12)$$

Therefore, first of all, we analyze the eigenvalues and eigenfunction of the operator $L_t(C)$. We assume that $C$ is the Hermitian matrix. Then the expression in (12) is the self-adjoint expression. Since the boundary conditions (2) are self-adjoint, the operator $L_t(C)$ is also self-adjoint. The eigenvalues of $C$, counted with multiplicity, and the corresponding orthonormal eigenvectors are denoted by $\mu_1 \leq \mu_2 \leq ... \leq \mu_m$ and $v_1, v_2, ..., v_m$. Thus

$$Cv_j = \mu_j v_j, \langle v_i, v_j \rangle = \delta_{i,j}.$$

One can easily verify that the eigenvalues and eigenfunctions of $L_t(C)$ are
$\mu_{k,j}(t) = (2\pi k + t)^n + \mu_j(2\pi k + t)^{n-2}$, $\Phi_{k,j,t}(x) = v_j e^{i(2\pi k + t)x}$, that is,

$$(L(C) - \mu_{k,j}(t))\Phi_{k,j,t}(x) = 0. \quad (13)$$

To prove the asymptotic formulas for the eigenvalues $\lambda_{k,j}(t)$ and for the corresponding normalized eigenfunctions $\Psi_{k,j,t}$ of $L_t(P)$ we use the formula

$$(\lambda_{k,j}(t) - \mu_{p,s}(t))(\Psi_{k,j,t}, \Phi_{p,s,t}) = (-i)^{n-2}((P_2 - C)\Psi_{k,j,t}^{(n-2)}, \Phi_{p,s,t}) + \sum_{\nu=3}^{n}(P_\nu \Psi_{k,j,t}^{n-\nu}, \Phi_{p,s,t})$$
$$(14)$$

which can be obtained from

$$L(P)\Psi_{k,j,t}(x) = \lambda_{k,j}(t)\Psi_{k,j,t}(x) \tag{15}$$

by multiplying both sides by $\Phi_{p,s,t}(x)$ and using (13). Then we estimate the right-hand side of (14) ( see Lemma 2) by using Lemma 1. At last, estimating $(\Psi_{k,j,t}, \Phi_{p,s,t})$ ( see Lemma 3) and using these estimations in (14), we find the asymptotic formulas for the eigenvalues and eigenfunctions of $L_t(P)$ ( see Theorem 1 and Theorem 2). Then using these formulas, we find the conditions on the eigenvalues of the matrix $C$ for which the number of the gaps in the spectrum of the operator $L(P)$ is finite ( see Theorem 3). Some of these results for differentiable $P_2(x)$ is obtained in [3, 6] by using the classical asymptotic expansions for the solutions of (6). The case $n = 2$ is investigated in [7]. In this paper we consider the more complicated case $n > 2$.

To estimate the right-hand side of (14) we use (8), (9), the following lemma, and the formula

$$\left(\lambda_{k,j}(t) - (2\pi p + t)^n\right)(\Psi_{k,j,t}, \varphi_{p,s,t}) = (-i)^{n-2}(P_2 \Psi_{k,j,t}^{n-2}, \varphi_{p,s,t}) + \sum_{\nu=3}^{n}(P_\nu \Psi_{k,j,t}^{n-\nu}, \varphi_{p,s,t}) \tag{16}$$

which can be obtained from (15) by multiplying both sides by $\varphi_{p,s,t}$ and using
$L_t(0)\varphi_{p,s,t} = (2\pi p + t)^n \varphi_{p,s,t}$.

**Lemma 1** *Let $\Psi_{k,j,t}(x)$ be normalized eigenfunction of $L_t(P)$. Then*

$$\sup_{x\in[0,1]} \left|\Psi_{k,j,t}^{(\nu)}(x)\right| = O(k^\nu) \tag{17}$$

*for $\nu = 0, 1, ..., n-2$. The equality (17) is uniform with respect to $t$ in $[-\frac{\pi}{2}, \frac{3\pi}{2})$.*

**Proof.** To prove (17) we use the arguments of the proof of the asymptotic formulas (4) and take into consideration the uniformity with respect to $t$. The eigenfunction $\Psi_{k,j,t}$ corresponding to the eigenvalue $\lambda_{k,j}(t)$ has the form

$$\Psi_{k,j,t}(x) = Y_1(x, \rho_{k,j})a_1 + Y_2(x, \rho_{k,j})a_2 + ... + Y_n(x, \rho_{k,j})a_n, \tag{18}$$

where $a_k \in \mathbb{C}^m$, $\rho_{k,j}(t) = i(\lambda_{k,j}(t))^{\frac{1}{n}}$, $Y_s(x, \rho_{k,j}(t))$ for $s = 1, 2, ..., n$ are linearly independent $m \times m$ matrix solutions of (6) for $\lambda = \lambda_{k,j}(t)$ satisfying

$$\frac{d^\nu Y_s(x, \rho_{k,j}(t))}{dx^\nu} = (\rho_{k,j}(t))^\nu e^{\rho_{k,j}(t)\omega_s x}\left[\omega_s^\nu I + O\left(\frac{1}{k}\right)\right] \tag{19}$$

for $\nu = 0, 1, ..., (n-1)$. Here $I$ is unit matrix, $\omega_1, \omega_2, ..., \omega_n$ are the $n$th root of 1, and $O\left(\frac{1}{k}\right)$ is a $m \times m$ matrix satisfying the following conditions

$$O\left(\frac{1}{k}\right) = \frac{A(x,t,k)}{k}, \quad |A(x,t,k)| < \frac{c_3}{|k|}, \quad \forall x \in [0,1], \forall t \in [-\frac{\pi}{2}, \frac{3\pi}{2}), \tag{20}$$

where $k > N$ and $c_3$ is a positive constant, independent of $t$. To consider the uniformity, with respect to $t$, of (17) we use (20).

*The proof of (17) in the case $n = 2r - 1$, $r > 1$.* Denote by $(\lambda_{k,j}(t))^{\frac{1}{n}}$ the root of $\lambda_{k,j}(t)$ lying in $O(k^{-\frac{1}{2m}})$ neighborhood of $(2k\pi + t)$ and put $\rho_{k,j}(t) = i(\lambda_{k,j}(t))^{\frac{1}{n}}$. Then we have

$$\rho_{k,j}(t) = (2k\pi + t)i + O(k^{-\frac{1}{2m}}). \tag{21}$$

Suppose $\omega_1, \omega_2, \ldots, \omega_n$ are ordered in such a way that

$$\omega_r = 1, \ \mathcal{R}(\rho_{k,j}(t)\omega_s) < 0, \forall s < r; \ \mathcal{R}(\rho_{k,j}(t)\omega_s) > 0, \forall s > r, \tag{22}$$

where $\mathcal{R}(z)$ is the real part of $z$. Using (19), (22), (2), (21), we get

$$Y_s^{(\nu-1)}(1, \rho_{k,j}(t)) = (\rho_{k,j}(t))^{\nu-1} e^{\rho_{k,j}(t)\omega_s}[\omega_s^{\nu-1}I], \ Y_s^{(\nu-1)}(0, \rho_{k,j}(t)) = (\rho_{k,j}(t))^{\nu-1}[\omega_s^{\nu-1}I],$$

$$U_\nu(Y_s(x, \rho_{k,j}(t))) = -(\rho_{k,j}(t))^{\nu-1} e^{it}[\omega_s^{\nu-1}I], \forall s < r,$$
$$U_\nu(Y_s(x, \rho_{k,j}(t))) = (\rho_{k,j}(t))^{\nu-1} e^{\rho_{k,j}(t)\omega_s}[\omega_s^{\nu-1}I], \forall s > r, \tag{23}$$

$$U_\nu(Y_r(x, \rho_{k,j}(t))) = (\rho_{k,j}(t))^{\nu-1} O(k^{-\frac{1}{2m}}), \tag{24}$$

where $[\omega_s^{\nu-1}I] = \omega_s^{\nu-1}I + O(\frac{1}{k})$ and $O(\frac{1}{k})$ satisfies the relation (20). Now using these relations and the notations of (18), we prove that

$$Y_s(x, \rho_{k,j}(t))a_s = O((\mid a_r \mid)k^{-\frac{1}{2m}}), \forall s \neq r. \tag{25}$$

Since $\Psi_{k,j,t}(x)$ satisfies (2) and (18), we have the system of equations

$$\sum_{s \neq r} U_\nu(Y_s(x, \rho_{k,j}(t)))a_s = -U_\nu(Y_r(x, \rho_{k,j}(t)))a_r, \ \nu = 0, 1, \ldots(n-2) \tag{26}$$

with respect to $a_{s,q}$ for $s \neq r$ and $q = 1, 2, \ldots, m$, where $a_{s,q}$ are coordinates of the vector $a_s$. Using (23) and (24) in (26) and then dividing both parts of $(\nu+1)$th equation of (26), for $\nu = 0, 1, \ldots(n-2)$, by $(\rho_{k,j}(t))^\nu$, we get the system of equations whose coefficient matrix $A$ is

$$\begin{pmatrix} -e^{it}[I] & \ldots & -e^{it}[I] & e^{\rho_{k,j}\omega_{r+1}}[I] & \ldots & e^{\rho_{k,j}\omega_n}[I] \\ -e^{it}[\omega_1 I] & \ldots & -e^{it}[\omega_{r-1} I] & e^{\rho_{k,j}\omega_{r+1}}[\omega_{r+1} I] & \ldots & e^{\rho_{k,j}\omega_n}[\omega_n I] \\ \ldots & \ldots & \ldots & \ldots & \ldots & \ldots \\ -e^{it}[\omega_1^{n-2} I] & \ldots & -e^{it}[\omega_{r-1}^{n-2} I] & e^{\rho_{k,j}\omega_{r+1}}[\omega_{r+1}^{n-2} I] & \ldots & e^{\rho_{k,j}\omega n}[\omega_n^{n-2} I] \end{pmatrix}$$

and the right-hand side is $O((\mid a_r \mid)k^{-\frac{1}{2m}})$. To estimate $\det A$ let us denote by $\widetilde{A}(m)$ the matrix obtained from $A$ by replacing $[\omega_s^j I]$ with $\omega_s^j I$ and by dividing the $s$th column ( note that the entries of the $s$th column are the $m \times m$ matrices) for $s < r$ and for $s > r$ by $-e^{it}$ and by $e^{\rho_{k,j}\omega_s}$ respectively. Clearly,

$$\widetilde{A}(1) = \begin{pmatrix} 1 & \ldots & 1 & 1 & \ldots & 1 \\ \omega_1 & \ldots & \omega_{r-1} & \omega_{r+1} & \ldots & \omega_n \\ \ldots & \ldots & \ldots & \ldots & \ldots & \ldots \\ \omega_1^{n-2} & \ldots & \omega_{r-1}^{n-2} & \omega_{r+1}^{n-2} & \ldots & \omega_n^{n-2} \end{pmatrix}$$

and $\det \widetilde{A}(1) \neq 0$. Besides, interchanging the rows and then interchanging the columns of $\widetilde{A}(m)$, we obtain $\det \widetilde{A}(m) = (\widetilde{A}(1))^m$. Using this and solving (26) by Cramer's rule, we get

$$a_{s,q} = \frac{\det A_{s,q}}{\det A} = O((\mid a_r \mid)e^{-\rho_{k,j}\omega_s}k^{-\frac{1}{2m}}), \forall s > r, \tag{27}$$

since $A_{s,q}$ is obtained from $A$ by replacing the $((s-1)m+q)$th column of $A$, which is the $q$th column of

$$\begin{pmatrix} e^{\rho_{k,j}\omega_s}[I] \\ e^{\rho_{k,j}\omega_s}[\omega_k I] \\ \ldots\ldots \\ e^{\rho_{k,j}\omega_s}[\omega_k^{n-2}] \end{pmatrix}, \text{ with } \begin{pmatrix} O(\mid a_r \mid k^{-\frac{1}{2m}}) \\ O(\mid a_r \mid k^{-\frac{1}{2m}}) \\ \ldots\ldots\ldots \\ O(\mid a_s \mid k^{-\frac{1}{2m}}) \end{pmatrix}. \text{ In the same way, we obtain}$$

$$a_{s,q} = O((\mid a_r \mid)k^{-\frac{1}{2m}}), \forall s < r. \tag{28}$$

Now (25) follows from (28), (27), (22). Therefore, the normalization condition $\parallel \Psi_{k,j,t} \parallel = 1$, and (25), (18), (19), (21), (22) imply that

$$\Psi_{k,j,t}(x) = Y_r(x, \rho_{k,j}(t)))a_r + O(k^{-\frac{1}{2m}}) = e^{i(2k\pi+t)x}a_r + O(k^{-\frac{1}{2m}}), \tag{29}$$

where $\mid a_r \mid^2 = 1 + O(k^{-\frac{1}{2m}})$, from which we get the proof of (17) for $\nu = 0$. Differentiating both sides of (18) and using (27), (28), we get the proof of (17) for arbitrary $\nu$ in the case $n = 2r - 1$.

*The proof of (17) in the case $n = 2r$.* In this case the $n$th roots $\omega_1, \omega_2, \ldots, \omega_n$ of 1 are ordered in such a way that

$$\omega_r = 1, \omega_{r+1} = -1, \mathcal{R}(\rho_{k,j}\omega_s) < 0, \forall s < r;\ \mathcal{R}(\rho_{k,j}\omega_s) > 0, \forall s > r+1. \tag{30}$$

Hence we have
$$U_\nu(Y_s(x, \rho_{k,j}(t))) = -(\rho_{k,j}(t))^{\nu-1}e^{it}[\omega_s^{v-1}I], \forall s < r, \tag{31}$$
$$U_\nu(Y_s(x, \rho_{k,j}(t))) = (\rho_{k,j}(t))^{\nu-1}e^{\rho_{k,j}(t)\omega_s}[\omega_s^{v-1}I], \forall s > r+1. \tag{32}$$

Now using these equalities, we prove that

$$Y_s(x, \rho_{k,j}(t))a_s = O((\mid a_r \mid + \mid a_{r+1} \mid)k^{-\frac{1}{2m}}), \forall s \neq r, r+1.$$

Using (31), (32) and arguing as in the case $n = 2r - 1$,, we get the system of equations

$$\sum_{s \neq r, r+1} U_\nu(Y_s(x, \rho_{k,j}(t)))a_s = -\sum_{s=r,r+1} U_\nu(Y_s(x, \rho_{k,j}(t)))a_s$$

for $\nu = 0, 1, 2, \ldots, (n-3)$. Arguing as in the proof of (27)-(29) and using (30), we get

$$a_{s,q} = O((\mid a_r \mid + \mid a_{r+1} \mid)e^{-\rho_{k,j}\omega_s}k^{-\frac{1}{2m}}), \forall s > r+1,$$

$$a_{s,q} = O((\mid a_r \mid + \mid a_{r+1} \mid)k^{-\frac{1}{2m}}), \forall s < r,$$

$$\Psi_{k,j,t}(x) = e^{i(2k\pi+t)x}a_r + e^{-i(2k\pi+t)x}a_{r+1} + O(k^{-\frac{1}{2m}}),$$

where $\mid a_r \mid^2 + \mid a_{r+1} \mid^2 = 1 + O(k^{-\frac{1}{2m}})$, which implies the proof of (17) in the case $n = 2r$
∎

It follows from this lemma that the equalities

$$(P_\nu \Psi_{k,j,t}^{n-\nu}, \varphi_{p,s,t}) = O(k^{n-\nu}),\ (P_\nu \Psi_{k,j,t}^{n-\nu}, \Phi_{p,s,t}) = O(k^{n-\nu}) \tag{33}$$

for $\nu = 2, 3, \ldots, n$ and for $j = 1, 2, \ldots, m$ hold uniformly with respect to $t$ in $[-\frac{\pi}{2}, \frac{3\pi}{2})$. Now (33) together with (16) implies that

$$\mid (\Psi_{k,j,t}, \varphi_{p,s,t}) \mid \leq \frac{c_4 \mid k \mid^{n-2}}{\mid \lambda_{k,j}(t) - (2\pi p + t)^n \mid} \tag{34}$$

for $p \notin A(k,n,t)$, $|k| \geq N$, and $s, j = 1, 2, ..., m$, where $c_4$ is a positive constant, independent of $t$. Using this we prove the following lemma.

**Lemma 2** *Let $b_{s,q}(x)$ be the entries of $P_2(x)$ and $b_{s,q,p} = \int_0^1 b_{s,q}(x) e^{-2\pi i p x} dx$. Then*

$$\left(\Psi_{k,j,t}^{(n-2)}, P_2 \varphi_{p,s,t}\right) = \sum_{q=1,2,...m;\ l \in \mathbb{Z}} b_{s,q,p-l}(\Psi_{k,j,t}^{(n-2)}, \varphi_{l,q,t}), \tag{35}$$

$$\left(\Psi_{k,j,t}^{(n-2)}, (P_2 - C)\Phi_{p,s,t}\right) = O(k^{n-3} \ln|k|) + O(k^{n-2} b_k) \tag{36}$$

*for $p \in A(k,n,t)$ and $s = 1, 2, ..., m$, where*

$$b_k = \max\{|b_{i,j,p}|: i, j = 1, 2, ...m;\ p = 2k, -2k, 2k+1, -2k-1\}, \tag{37}$$

*and $C$ is the Hermitian matrix defined in (12). The formula (36) is uniform with respect to $t$ in $[-\frac{\pi}{2}, \frac{3\pi}{2})$. Moreover, in the Case 1 of Notation 1 the formula*

$$\left(\Psi_{k,j,t}^{(n-2)}, (P_2 - C)\Phi_{k,s,t}\right) = O(k^{n-3} \ln|k|) \tag{38}$$

*holds. If $n = 2r - 1$, then (38) is uniform with respect to $t$ in $[-\frac{\pi}{2}, \frac{3\pi}{2})$.*

**Proof.** Note that if the entries $b_{s,q}$ of $P_2$ belong to $L_2[0,1]$, then (35) is obvious, since $\{\varphi_{l,q,t}:\ l \in \mathbb{Z},\ q = 1, 2, ..., m\}$ is an orthonormal basis in $L_2^m[0,1]$. Now we prove (35) in case $b_{s,q} \in L_1[0,1]$. Using (2), (34) and the integration by parts, we see that there exists a constant $c_5$, independent of $t$, such that

$$\left|\left(\Psi_{k,j,t}^{(n-2)}, \varphi_{l,q,t}\right)\right| = |(2\pi l + t)^{n-2} (\Psi_{k,j,t}, \varphi_{l,q,t})| \leq \frac{c_5 |k|^{n-2} |l|^{n-2}}{|\lambda_{k,j}(t) - (2\pi l + t)^n|}. \tag{39}$$

for $l \notin A(k,n,t)$, $|k| \geq N$. This and (8) imply that there exists a constant $c_6$, independent of $t$, such that

$$\sum_{l: l > d} \left|\left(\Psi_{k,j,t}^{(n-2)}, \varphi_{l,q,t}\right)\right| < \frac{c_6 |k|^{n-2}}{d},$$

where $d > 2|k|$, $t \in [-\frac{\pi}{2}, \frac{3\pi}{2})$. Hence the decomposition of $\Psi_{k,j,t}^{(n-2)}$ by the basis $\{\varphi_{l,q,t}:\ l \in \mathbb{Z},\ q = 1, 2, ..., m\}$ has the form

$$\Psi_{k,j,t}^{(n-2)}(x) = \sum_{|l| \leq d;\ q=1,2,...,m} \left(\Psi_{k,j,t}^{(n-2)}, \varphi_{l,q,t}\right) \varphi_{l,q,t}(x) + g_d(x), \tag{40}$$

$$\text{where} \quad \sup_{x \in [0,1]} |g_d(x)| < \frac{c_6 |k|^{n-2}}{d}.$$

Using (40) in $\left(\Psi_{k,j,t}^{(n-2)}, P_2 \varphi_{p,s,t}\right)$ and letting $d$ tend to $\infty$, we obtain (35).

Since $\Phi_{p,s,t}(x) \equiv v_s e^{i(2\pi p + t)x}$, to prove (36), it is enough to show that

$$\left(\Psi_{k,j,t}^{(n-2)}, (P_2 - C)\varphi_{p,s,t}\right) = O((k^{n-3} \ln|k|) + O(k^{n-2} b_k)$$

for $s = 1, 2, ..., m$ and $p \in A(k,n,t)$. Using the obvious relation

$$\left(\Psi_{k,j,t}^{(n-2)}, C\varphi_{p,s,t}\right) = \sum_{q=1,2,...,m} b_{s,q,0}(\Psi_{k,j,t}^{(n-2)}, \varphi_{p,q,t})$$

and (35), we see that

$$\left(\Psi_{k,j,t}^{(n-2)}, (P_2 - C)\varphi_{p,s,t}\right) = \sum_{l: l \in A(k,n,t)\setminus p;\ q=1,2,\ldots,m} b_{s,q,p-l} \left(\Psi_{k,j,t}^{(n-2)}, \varphi_{l,q,t}\right) \quad (41)$$

$$+ \sum_{l: l \notin A(k,n,t);\ q=1,2,\ldots,m} b_{s,q,p-l} \left(\Psi_{k,j,t}^{(n-2)}, \varphi_{l,q,t}\right).$$

Since

$$\mid b_{j,i,s} \mid \leq \max_{p,q=1,2,\ldots,m} \int_0^1 \mid b_{p,q}(x) \mid dx = O(1)$$

for all $j$, $i$, $s$, using (39) and (9), we see that the second summation of the right-hand side of (41) is $O((k^{n-3}\ln|k|)$. Besides, it follows from (17), (37) that the first summation of the right-hand side of (41) is $O(k^{n-2}b_k)$, since for $p \in A(k,n,t)$ and $l \in A(k,n,t)\setminus p$, we have $p - l \in \{2k, -2k, 2k+1, -2k-1\}$. Hence (36) is proved. In the Case 1 of Notation 1 the first summation of the right-hand side of (41) is absent, since in this case $A(k,n,t) = \{k\}$ and $A(k,n,t)\setminus p = \emptyset$ for $p \in A(k,n,t)$. Thus (38) is proved. The uniformity of the formulas (36), (38) follows from the uniformity of (17), (8) and (9). ∎

**Lemma 3** *There exists a positive number $N_0$, independent of $t$, such that for $\mid k \mid > N_0$ and for $p \in A(k,n,t)$ the following assertions hold.*

*(a) If $C$ is Hermitian matrix, then for each eigenfunction $\Psi_{k,j,t}$ of $L_t(P)$ there exists an eigenfunction $\Phi_{p,s,t}$ of $L_t(C)$ satisfying*

$$|(\Psi_{k,j,t}, \Phi_{p,s,t})| > \frac{1}{3m}. \quad (42)$$

*(b) If $L_t(P)$ is self adjoint operator, then for each eigenfunction $\Phi_{k,j,t}$ of $L_t(C)$ there exists an eigenfunction $\Psi_{p,s,t}$ of $L_t(P)$ satisfying*

$$|(\Phi_{k,j,t}, \Psi_{p,s,t})| > \frac{1}{3m}. \quad (43)$$

**Proof.** It follows from (34) and (10) that

$$\sum_{s=1,2,\ldots,m} \left( \sum_{p: p \notin A(k,n,t)} |(\Psi_{k,j,t}, \varphi_{p,s,t})|^2 \right) = O\left(\frac{(\ln|k|)^2}{k^2}\right).$$

Hence using the equality $\Phi_{p,s,t}(x) = v_s e^{i(2\pi p + t)x}$, where $v_s$ is the normalized eigenvectors of $C$, and the Parseval equality, we get

$$\sum_{s=1,2,\ldots,m} \left( \sum_{p: p \notin A(k,n,t)} |(\Psi_{k,j,t}, \Phi_{p,s,t})|^2 \right) = O\left(\frac{(\ln|k|)^2}{k^2}\right), \quad (44)$$

$$\sum_{s=1,2,\ldots,m;\ p \in A(k,n,t)} |(\Psi_{k,j,t}, \Phi_{p,s,t})|^2 = 1 + O\left(\frac{(\ln|k|)^2}{k^2}\right). \quad (45)$$

Since the number of the eigenfunctions $\Phi_{p,s,t}(x)$ for $p \in A(k,n,t)$, $s = 1, 2, \ldots, m$ is less than $2m$ (see Notation 1), (42) follows from (45).

Using (34) and (11), we get

$$\sum_{s=1,2,...,m}(\sum_{p:p\notin A(k,n,t)} |(\varphi_{k,j,t}, \Psi_{p,s,t})|^2) = O(\frac{(\ln|k|)^2}{k^2})$$

Therefore, arguing as in the proof of (42) and taking into account that the eigenfunctions of the self-adjoint operator $L_t(P)$ form an orthonormal basis in $L_2^m(0,1)$, we get the proof of (43) ∎

**Theorem 1** *Let $L_t(P)$ be a self-adjoint operator and $C$ be a Hermitian matrix. If $n = 2r-1$, then for arbitrary $t$, if $n = 2r$, then for $t \neq 0, \pi$ the large eigenvalues of $L_t(P)$ consist of $m$ sequences (3) satisfying*

$$\lambda_{k,j}(t) = (2\pi k + t)^n + \mu_j (2\pi k + t)^{n-2} + O(k^{n-3}\ln|k|) \tag{46}$$

*and the normalized eigenfunction $\Psi_{k,j,t}$ corresponding to $\lambda_{k,j}(t)$ satisfies*

$$\| \Psi_{k,j,t} - E\Psi_{k,j,t} \| = O(\frac{(\ln|k|)}{k}) \tag{47}$$

*for $j = 1, 2, ..., m$, where $\mu_1 \leq \mu_2 \leq ... \leq \mu_m$ are the eigenvalues of $C$ and $E$ is the orthogonal projection onto the eigenspace of $L_t(C)$ corresponding to $\mu_{k,j}(t)$. If $\mu_j$ is a simple eigenvalue of $C$, then the eigenvalue $\lambda_{k,j}(t)$ satisfying (46) is a simple eigenvalue, and the corresponding eigenfunction satisfies*

$$\Psi_{k,j,t}(x) = v_j e^{i(2\pi k+t)x} + O(\frac{(\ln|k|)}{k}), \tag{48}$$

*where $v_j$ is the eigenvector of $C$ corresponding to the eigenvalue $\mu_j$. In the case $n = 2r - 1$ the formulas (46)-(48) are uniform with respect to $t$ in $[-\frac{\pi}{2}, \frac{3\pi}{2})$.*

**Proof.** By (33) and (38) the right-hand side of (14) is $O(k^{n-3}\ln|k|)$. On the other hand by Notation 1 if $t \neq 0, \pi$, then there exists $N$ such that $t \in T(k)$, and hence $A(k,n,t) = \{k\}$, for $|k| \geq N$. Thus dividing (14) by $(\Psi_{k,j,t}, \Phi_{p,s,t})$, where $p \in A(k,n,t)$, and hence $p = k$, and using (42), we get

$$\{\lambda_{k,1}(t), \lambda_{k,2}(t), ..., \lambda_{k,m}(t)\} \subset \cup_{j=1}^m (U(\mu_{k,j}(t), \delta_k)), \tag{49}$$

where $U(\mu, \delta) = \{z \in \mathbb{R} : |\mu - z| < \delta\}$, $|k| \geq \max\{N, N_0\}$, $\delta_k = O(|k|^{n-3}\ln|k|)$. Instead of (42) using (43), in the same way, we obtain

$$U(\mu_{k,s}(t), \delta_k) \cap \{\lambda_{k,1}(t), \lambda_{k,2}(t), ..., \lambda_{k,m}(t)\} \neq \emptyset \tag{50}$$

for $|k| \geq \max\{N, N_0\}$ and $s = 1, 2, ..., m$. Hence to prove (46) we need to show that if the multiplicity of the eigenvalue $\mu_j$ is $q$ then there exist precisely $q$ eigenvalues of $L_t(P)$ lying in $U(\mu_{k,j}(t), \delta_k)$ for $|k| \geq \max\{N, N_0\}$. The eigenvalues of $L_t(P)$ and $L_t(C)$ can be numbered in the following way: $\lambda_{k,1}(t) \leq \lambda_{k,2}(t) \leq ... \leq \lambda_{k,m}(t)$ and $\mu_{k,1}(t) \leq \mu_{k,2}(t) \leq ... \leq \mu_{k,m}(t)$. If $C$ has $\nu$ different eigenvalues $\mu_{j_1}, \mu_{j_2}, ..., \mu_{j_\nu}$ with multiplicities $j_1, j_2 - j_1, ..., j_\nu - j_{\nu-1}$, then we have

$$j_1 < j_2 < ... < j_\nu = m;\ \mu_{j_1} < \mu_{j_2} < ... < \mu_{j_\nu};\ \mu_1 = \mu_2 = ... = \mu_{j_1}; \tag{51}$$
$$\mu_{j_1+1} = \mu_{j_1+2} = ... = \mu_{j_2}; ...; \mu_{j_{\nu-1}+1} = \mu_{j_{\nu-2}+2} = ... = \mu_{j_\nu}.$$

Suppose there exist precisely $s_1, s_2, ..., s_\nu$ eigenvalues of $L_t(P)$ lying in the intervals

$$U(\mu_{k,j_1}(t), \delta_k), U(\mu_{k,j_2}(t), \delta_k), ..., U(\mu_{k,j_\nu}(t), \delta_k)$$

respectively. Since

$$\delta_k \ll (\min_{p=1,2,...,\nu-1} \mid (\mu_{j_p+1} - \mu_{j_p})(2\pi k + t)^{n-2} \mid) \text{ for } \mid k \mid \gg 1,$$

these intervals are pairwise disjoints. Therefore using (3) and (4), we get

$$s_1 + s_2 + ... + s_\nu = m. \tag{52}$$

Now let us prove that $s_1 = j_1$, $s_2 = j_2 - j_1$, ..., $s_\nu = j_\nu - j_{\nu-1}$. Due to the notations the eigenvalues $\lambda_{k,1}(t), \lambda_{k,2}(t), ..., \lambda_{k,s_1}(t)$ of the operator $L_t(P)$ lie in $U(\mu_{k,1}(t), \delta_k)$ and by the definition of $\delta_k$ we have

$$\mid \lambda_{k,j}(t) - \mu_{k,s}(t) \mid > \frac{1}{2}(\min_{p:p>j_1} \mid (\mu_1 - \mu_p)(2\pi k + t)^{n-2} \mid)$$

for $j \leq s_1$ and $s > j_1$. Hence using (14) for $p = k$ and (38), (33), we get

$$\sum_{s: s > j_1} \mid (\Psi_{k,j,t}, \Phi_{k,s,t}) \mid^2) = O(\frac{(\ln \mid k \mid)^2}{k^2}), \ \forall j \leq s_1.$$

Using this, (44), and taking into account that $A(k, n, t) = \{k\}$ for $\mid k \mid \geq N$, we conclude that there exists normalized eigenfunction, denoted by $\Phi_{k,j,t}(x)$, of $L_t(P)$ corresponding to $\mu_{k,1}(t) = \mu_{k,2}(t) = ... = \mu_{k,j_1}(t)$ such that

$$\Psi_{k,j,t}(x) = \Phi_{k,j,t}(x) + O(k^{-1} \ln \mid k \mid) \tag{53}$$

for $j \leq s_1$. Since $\Psi_{k,1,t}, \Psi_{k,2,t}, ..., \Psi_{k,s_1,t}$ are orthonormal system we have

$$(\Phi_{k,j,t}, \Phi_{k,s,t}) = \delta_{s,j} + O(k^{-1} \ln \mid k \mid), \ \forall s, j = 1, 2, ..., s_1.$$

This formula implies that the dimension $j_1$ of the eigenspace of $L_t(C)$ corresponding to the eigenvalue $\mu_{k,1}(t)$ is not less than $s_1$. Thus $s_1 \leq j_1$. In the same way we prove that

$s_2 \leq j_2 - j_1, ..., s_\nu \leq j_\nu - j_{\nu-1}$. Now (52) and the equality $j_\nu = m$ (see (51)) imply that

$s_1 = j_1, s_2 = j_2 - j_1, ..., s_\nu = j_\nu - j_{\nu-1}$. Therefore, taking into account that, the eigenvalues of $L_t(P)$ consist of $m$ sequences satisfying (4), we get (46). The proof of (47) follows from (53).

Now suppose that $\mu_j$ is a simple eigenvalue of $C$. Then $\mu_{k,j}(t)$ is a simple eigenvalues of $L_t(C)$ and, as it was proved above, there exists unique eigenvalues $\lambda_{k,j}(t)$ of $L_t(P)$ lying in $U(\mu_{k,s}(t), \delta_k)$, where $\mid k \mid \geq \max\{N, N_0\}$, and the eigenvalues $\lambda_{k,j}(t)$ for

$\mid k \mid \geq \max\{N, N_0\}$ are the simple eigenvalues. Hence (48) is the consequence of (47), since there exists unique eigenfunction $\Phi_{k,j,t}(x) = v_j e^{i(2\pi k + t)x}$ corresponding to the eigenvalue $\mu_{k,j}(t)$. The uniformity of the formulas (46)-(48) follows from the uniformity of (38), (33), (42), (43) ∎

**Theorem 2** *Let $L_t(P)$ be a self-adjoint operator, $C$ be a Hermitian matrix, $n = 2r$, $\mu_j$ be a simple eigenvalue of $C$, $\alpha_j$ be a positive constant satisfying $\alpha_j < \min_{q:q \neq j} \mid \mu_j - \mu_q \mid$ and* $B(\alpha_j, k, \mu_j)$ be a set defined by $B(\alpha_j, k, \mu_j) = B(0, \alpha_j, k, \mu_j) \cup B(\pi, \alpha_j, k, \mu_j)$, where

$$B(0, \alpha_j, k, \mu_j) = \bigcup_{s=1,2,\ldots,m} (\frac{\mu_s - \mu_j - \alpha_j}{4n\pi k}, \frac{\mu_s - \mu_j + \alpha_j}{4n\pi k}),$$

$$B(\pi, \alpha_j, k, \mu_j) = \bigcup_{s=1,2,\ldots,m} (\pi + \frac{\mu_s - \mu_j - \alpha_j}{2n\pi(2k+n-1)}, \pi + \frac{\mu_s - \mu_j + \alpha_j}{2n\pi(2k+n-1)}).$$

*There exist a positive number $N_1$ such that if $\mid k \mid \geq N_1$ and $t \notin B(\alpha_j, k, \mu_j)$, then there exists a unique eigenvalues, denoted by $\lambda_{k,j}(t)$, of $L_t(P)$ lying in $U(\mu_{k,j}, \varepsilon_k)$, where*
$\varepsilon_k = c_7(|k|^{n-3} \ln|k|) + |k|^{n-2} b_k$, $b_k$ *is defined by (37), and $c_7$ is a positive constant, independent of $t$. The eigenvalue $\lambda_{k,j}(t)$ is a simple eigenvalue of $L_t(P)$ and the corresponding normalized eigenfunction $\Psi_{k,j,t}(x)$ satisfies*

$$\Psi_{k,j,t}(x) = v_j e^{i(2\pi k + t)x} + O(k^{-1} \ln|k|) + O(b_k). \tag{54}$$

**Proof.** To consider the simplicity of $\mu_{k,j}(t)$ and $\lambda_{k,j}(t)$ we introduce the set

$$S(k,j,p,s) = \{t \in [-\frac{\pi}{2}, \frac{3\pi}{2}) : \mid \mu_{k,j}(t) - \mu_{p,s}(t) \mid < \alpha_j |k|^{n-2}\} \tag{55}$$

for $(p,s) \neq (k,j)$. It follows from (7) that $S(k,j,p,s) = \emptyset$ for $p \neq k, -k, -k-1$. Moreover, if $\mu_j$ is a simple eigenvalue, then $S(k,j,k,s) = \emptyset$ for $s \neq j$, since

$$\mid \mu_{k,j}(t) - \mu_{k,s}(t) \mid = \mid (\mu_j - \mu_s)(2\pi k + t)^{n-2} \mid > \alpha_j |k|^{n-2}.$$

It remains to consider the sets $S(k,j,-k,s), S(k,j,-k-1,s)$. Using the equality
$\mu_{k,j}(t) - \mu_{-k,s}(t) = (2\pi k)^{n-2}(4nk\pi t + \mu_j - \mu_s) + O(k^{n-3})$, we see that

$$S(k,j,-k,s) \subset (\frac{\mu_s - \mu_j - \alpha_j}{4n\pi k}, \frac{\mu_s - \mu_j + \alpha_j}{4n\pi k}).$$

Similarly, by using the obvious equality $\mu_{k,j}(t) - \mu_{-k-1,s}(t) =$
$(2\pi k + t)^n + \mu_j(2\pi k + t)^{n-2} - (2\pi k + 2\pi - t)^n - \mu_j(2\pi k + 2\pi - t)^{n-2} =$

$$(2\pi k)^{n-2}(n2\pi kt - n2\pi k(2\pi - t) + \frac{1}{2}n(n-1)(t^2 - (2\pi - t)^2) + \mu_j - \mu_s) + O(k^{n-3}) =$$

$(2\pi k)^{n-2}[(t-\pi)(2k+(n-1))2\pi n + \mu_j - \mu_s] + O(k^{n-3})$, we get

$$S(k,j,-k-1,s) \subset (\pi + \frac{\mu_s - \mu_j - \alpha_j}{2n\pi(2k+n-1)}, \pi + \frac{\mu_s - \mu_j + \alpha_j}{2n\pi(2k+n-1)}).$$

Using these relations and the definition of $B(\alpha_j, k, \mu_j)$, we obtain

$$\bigcup_{\substack{p \in \mathbb{Z}, s=1,2,\ldots,m, \\ (p,s) \neq (k,j)}} S(k,j,p,s) = \bigcup_{\substack{p=-k,-k-1 \\ s=1,2,\ldots,m}} S(k,j,p,s) \subset B(\alpha_j, k, \mu_j).$$

Therefore it follows from (55) that if $t \notin B(\alpha_j, k, \mu_j)$, then

$$\mid \mu_{k,j}(t) - \mu_{p,s}(t) \mid \geq \alpha_j |k|^{n-2} \tag{56}$$

for all $(p,s) \neq (k,j)$. Hence $\mu_{k,j}(t)$ is a simple eigenvalue of $L_t(C)$ for $t \notin B(\alpha_j, k, \mu_j)$. Instead of (38) using (36) and arguing as in the proof of (50), we obtain that there exists $N_1$ such that if $\mid k \mid \geq N_1$, then there exists an eigenvalue, denoted by $\lambda_{k,j}(t)$, of $L_t(P)$ lying in

$U(\mu_{k,j}(t), \varepsilon_k)$. Now using the definition of $\varepsilon_k$ and then (56), we see that

$$\mid \lambda_{k,j}(t) - \mu_{k,j}(t) \mid < \varepsilon_k = o(k^{n-2}), \quad \mid \lambda_{k,j}(t) - \mu_{p,s}(t) \mid > \frac{1}{2}\alpha_j|k|^{n-2}. \tag{57}$$

for $\mid k \mid \geq N_1$, $s = 1, 2, ..., m$, $(p, s) \neq (k, j)$ and for any eigenvalue $\lambda_{k,j}(t)$ lying in $U(\mu_{k,j}(t), \varepsilon_k)$. Let $\Psi_{k,j,t}(x)$ be any normalized eigenfunction corresponding to $\lambda_{k,j}(t)$. Dividing both sides of (14) by $\lambda_{k,j}(t) - \mu_{p,s}(t)$ and using (36), (33), (57), we get

$$(\Psi_{k,j,t}, \Phi_{p,s,t}) = O(k^{-1} \ln |k|) + O(b_k)$$

for $(p, s) \neq (k, j)$ and $p \in A(k, n, t)$. This and (44), (45) imply that $\Psi_{k,j,t}(x)$ satisfies (54). Thus we have proved that (54) holds for any normalized eigenfunction of $L_t(P)$ corresponding to any eigenvalue lying in $U(\mu_{k,j}(t), \varepsilon_k)$. If there exist two different eigenvalues of $L_t(P)$ lying in $U(\mu_{k,j}(t), \varepsilon_k)$ or if there exists a multiple eigenvalue of $L_t(P)$ lying in $U(\mu_{k,j}(t), \varepsilon_k)$, then we obtain that there exist two orthonormal eigenfunctions satisfying (54) which is impossible. Therefore there exists unique eigenvalue $\lambda_{k,j}(t)$ of $L_t(P)$ lying in $U(\mu_{k,j}(t), \varepsilon_k)$ and $\lambda_{k,j}(t)$ is a simple eigenvalue of $L_t(P)$ ∎

**Theorem 3** *Let $L(P)$ be self-adjoint operator generated in $L_2^m(-\infty, \infty)$ by the differential expression (1) and $C$ be Hermitian matrix.*

*(a) If $n$ and $m$ are odd numbers then the spectrum $\sigma(L(P))$ of $L(P)$ coincides with $(-\infty, \infty)$.*

*(b) If $n$ is odd number, $n > 1$, and the matrix $C$ has at least one simple eigenvalue, then the number of the gaps in $\sigma(L(P))$ is finite.*

*(c) Suppose that $n$ is even number, and the matrix $C$ has at least three simple eigenvalues $\mu_{j_1} < \mu_{j_2} < \mu_{j_3}$ such that $diam(\{\mu_{j_1} + \mu_{i_1}, \mu_{j_2} + \mu_{i_2}, \mu_{j_3} + \mu_{i_3}\}) \neq 0$ for each triple $(i_1, i_2, i_3)$, where $i_p = 1, 2, ..., m$ for $p = 1, 2, 3$ and $diam(A)$ is the diameter $\sup_{x,y \in A} \mid x - y \mid$ of the set $A$. Then the number of the gaps in the spectrum of $L(P)$ is finite.*

**Proof.** (a) In case $m = 1$ the assertion (a) is proved in [4]. Our proof is carried out analogous fashion. Since $L(P)$ is self-adjoint, $\sigma(L(P))$ is a subset of $(-\infty, \infty)$. Therefore we need to prove that $(-\infty, \infty) \subset \sigma(L(P))$. Suppose to the contrary that there exists a real number $\lambda$ such that $\lambda \notin \sigma(L(P))$. It is not hard to see that the characteristic determinant $\Delta(\lambda, t) = \det(U_\nu(Y_s(x, \lambda)))$ of $L_t(P)$ has the form

$$\Delta(\lambda, t) = e^{inmt} + a_1(\lambda)e^{i(nm-1)t} + a_2(\lambda)e^{i(nm-2)t} + ... + a_{nm}(\lambda), \tag{58}$$

that is, $\Delta(\lambda, t)$ is a polynomial $S_\lambda(u)$ of $u = e^{it}$ of order $nm$ with entire coefficients $a_1(\lambda)$, $a_2(\lambda)$, .... It is well known that if $\lambda \notin \sigma(L(P))$, then the absolute values of all roots $u_1 = e^{it_1}$, $u_2 = e^{it_2}, ..., u_{nm} = e^{it_{nm}}$ of $S_\lambda(u) = 0$ differ from 1, that is, $t_k \neq \overline{t_k}$ and $\lambda$ is the eigenvalue of $L_{t_k}(P)$ for $k = 1, 2, ..., nm$. It is not hard to see that $L_{t_k}^* = L_{\overline{t_k}}$, $\lambda = \overline{\lambda} \in \sigma(L_{\overline{t_k}})$. Moreover, if $\lambda$ is the eigenvalue of $L_{t_k}(P)$ of multiplicity $m_k$ then $\overline{\lambda}$ is the eigenvalue of $L_{\overline{t_k}}(P)$ of the same multiplicity $m_k$. Now taking into account that $u_k = e^{it_k}$ is the root of $S_\lambda(u) = 0$ of multiplicity $m_k$ if and only if $\lambda$ is the eigenvalue of $L_{t_k}(P)$ of multiplicity $m_k$, we obtain that $e^{i\overline{t_k}}$ is also root of $S_\lambda(u) = 0$ of the same multiplicity $m_k$. Since $e^{i\overline{t_k}} \neq e^{it_k}$, we see that the number $nm$ of the roots of $S_\lambda(u) = 0$ ( see (58)) is an even number which contradicts the assumption that $n$ and $m$ are odd numbers.

(b) It follows from the uniform asymptotic formula (46) that there exists a positive numbers $N_2$, $c_8$, independent of $t$, such that if $\mid k \mid \geq N_2$ and $\mu_j$ is a simple eigenvalue of the matrix $C$ then there exists unique simple eigenvalue $\lambda_{k,j}(t)$ of $L_t(P)$ lying in $U(\mu_{k,j}(t), \delta_k)$, where $\delta_k = c_8|k|^{n-3} \ln |k|$ and $t \in [-\frac{\pi}{2}, \frac{3\pi}{2})$. Therefore $\lambda_{k,j}(t_0)$ for $t_0 \in (-\frac{\pi}{2}, \frac{3\pi}{2})$, $\mid k \mid \geq N_2$ is a simple zero of the characteristic determinant $\Delta(\lambda, t_0)$. By implicit function theorem there exists a neighborhood $U(t_0) \subset (-\frac{\pi}{2}, \frac{3\pi}{2})$ of $t_0$ and a continuous in $U(t_0)$ function $\Lambda(t)$ such that $\Lambda(t_0) = \lambda_{k,j}(t_0)$, $\Lambda(t)$ is an eigenvalue of $L_t(P)$ for $t \in U(t_0)$ and $|\Lambda(t) - \mu_{k,j}(t)| < \delta_k$, for all $t \in U(t_0)$, since $|\Lambda(t_0) - \mu_{k,j}(t_0)| = |\lambda_{k,j}(t_0) - \mu_{k,j}(t_0)| < \delta_k$ and the functions $\Lambda(t), \mu_{k,j}(t)$ are continuous. Now taking into account that there exists unique eigenvalue of $L_t(P)$ lying in $U(\mu_{k,j}(t), \delta_k)$, we obtain that $\Lambda(t) = \lambda_{k,j}(t)$ for $t \in U(t_0)$, and hence $\lambda_{k,j}(t)$ is continuous at $t_0 \in (-\frac{\pi}{2}, \frac{3\pi}{2})$. Therefore the sets $\Gamma_{k,j} = \{\lambda_{k,j}(t) : t \in (-\frac{\pi}{2}, \frac{3\pi}{2})\}$ for $|k| \geq N_2$ are intervals and $\Gamma_{k,j} \subset \sigma(L(P))$. Similarly there exists a neighborhood

$$U(-\frac{\pi}{2}) = (-\frac{\pi}{2} - \beta, -\frac{\pi}{2} + \beta)$$

of $-\frac{\pi}{2}$ and a continuous function $A(t)$ such that $A(-\frac{\pi}{2}) = \lambda_{k+1,j}(-\frac{\pi}{2})$, where $k \geq N_2$, $A(t)$ is an eigenvalue of $L_{t+2\pi}(P)$ ($L_t(P) = L_{t+2\pi}(P)$) for $t \in (-\frac{\pi}{2} - \beta, -\frac{\pi}{2}]$ and $A(t)$ is an eigenvalue of $L_t(P)$ for $t \in [-\frac{\pi}{2}, -\frac{\pi}{2} + \beta)$ and

$$|A(t) - \mu_{k+1,j}(t)| < \delta_k \text{ for all } t \in [-\frac{\pi}{2}, -\frac{\pi}{2} + \beta),$$
$$|A(t) - \mu_{k,j}(t+2\pi)| < \delta_k \text{ for all } t \in (-\frac{\pi}{2} - \beta, -\frac{\pi}{2}), \text{ since}$$

$|A(-\frac{\pi}{2}) - \mu_{k+1,j}(-\frac{\pi}{2})| = |\lambda_{k+1,j}(-\frac{\pi}{2}) - \mu_{k+1,j}(-\frac{\pi}{2})| < \delta_k$, $\mu_{k,j}(t+2\pi) = \mu_{k+1,j}(t)$ and the functions $\Lambda(t), \mu_{k,j}(t)$ are continuous. Again taking into account that there exists unique eigenvalue of $L_t(P)$ lying in $U(\mu_{k+1,j}(t), \delta_k)$ for $t \in [-\frac{\pi}{2}, -\frac{\pi}{2} + \beta)$ and lying in $U(\mu_{k,j}(t), \delta_k)$ for $t \in (\frac{3\pi}{2} - \beta, \frac{3\pi}{2})$, we obtain that

$$A(t) = \lambda_{k+1,j}(t), \, \forall t \in [-\frac{\pi}{2}, -\frac{\pi}{2} + \beta) \text{ and } A(t) = \lambda_{k,j}(t+2\pi), \, \forall t \in (-\frac{\pi}{2} - \beta, -\frac{\pi}{2}).$$

Thus one part of the interval $\{A(t) : t \in (-\frac{\pi}{2} - \beta, -\frac{\pi}{2} + \beta)\}$ lies in $\Gamma_{k,j}$ and the other part lies in $\Gamma_{k+1,j}$, that is, the interval $\Gamma_{k,j}$ and $\Gamma_{k+1,j}$ are connected for $k \geq N_2$. Similarly the interval $\Gamma_{k,j}$ and $\Gamma_{k-1,j}$ are connected for $k \leq -N_2$. Therefore the number of the gaps in the spectrum of $L(P)$ is finite.

(c) In Theorem 2 we proved that if $|k| \geq N_1$ and $t \notin B(\alpha_{j_p}, k, \mu_{j_p})$, where $p = 1, 2, 3$, then there exists a unique eigenvalue, denoted by $\lambda_{k,j_p}(t)$, of $L_t(P)$ lying in $U(\mu_{k,j_p}(t), \varepsilon_k)$ and it is a simple eigenvalue. Let us prove that $\lambda_{k,j_p}(t)$ is continuous at

$$t_0 \in [-\frac{\pi}{2}, \frac{3\pi}{2}) \backslash B(\alpha_{j_p}, k, \mu_{j_p}).$$

Since $\lambda_{k,j_p}(t_0)$ is a simple eigenvalue it is a simple zero of the characteristic determinant $\Delta(\lambda, t)$ of the operator $L_t(P)$. Therefore repeating the argument of the proof of the continuity of $\lambda_{k,j}(t)$ in the proof of (b), we obtain that $\lambda_{k,j_p}(t)$ is continuous at $t_0$ for $|k| \geq N_1$. Now we prove that there exists $H$ such that

$$(H, \infty) \subset \{\lambda_{k,j_p}(t) : t \in [-\frac{\pi}{2}, \frac{3\pi}{2}) \backslash B(\alpha_{j_p}, k, \mu_{j_p}), \, k = N_1, N_1 + 1, ..., \}.$$

It is clear that

$$(h, \infty) \subset \{\mu_{k,j_p}(t) : t \in [-\frac{\pi}{2}, \frac{3\pi}{2}), k = N_1, N_1 + 1, ...\}, \, \forall p = 1, 2, 3, \quad (59)$$

where $h = \mu_{N_1, j_3}(-\frac{\pi}{2})$. Since $\mu_{k,j_p}(t)$ is increasing function for $k \geq N_1$, it follows from the

obvious equality

$$\mu_{k,j_p}(\frac{\mu_s - \mu_{j_p} \mp \alpha_{j_p}}{4n\pi k}) = (2\pi k)^n + n(2\pi k)^{n-1}\frac{\mu_s - \mu_{j_p} \mp \alpha_{j_p}}{4n\pi k} + \mu_{j_p}(2\pi k)^{n-2}$$

$$+O(k^{n-4}) = (2\pi k)^n + (2\pi k)^{n-2}\frac{\mu_s + \mu_{j_p} \mp \alpha_{j_p}}{2} + O(k^{n-4})$$

and from the definition of $B(0, \alpha_{j_p}, k, \mu_{j_p})$ that

$$\{\mu_{k,j_p}(t) : t \in B(0, \alpha_{j_p}, k, \mu_{j_p})\} \subset \bigcup_{s=1,2\ldots,m} C(0, k, j_p, s, \alpha_{j_p}), \text{ where}$$

$C(0, k, j_p, s, \alpha_{j_p}) = \{x \in \mathbb{R} : |x - (2\pi k)^n + (2\pi k)^{n-2}\frac{\mu_s + \mu_{j_p}}{2}| < \alpha_{j_p}(2\pi k)^{n-2}\}$. This inclusion with (59) implies that the set

$$(h, \infty) \backslash \bigcup_{k: k \geq N_1;\ s=1,2\ldots,m} C(0, k, j_p, s, \alpha_{j_p})$$

is a subset of the set $\{\mu_{k,j_p}(t) : t \in [-\frac{\pi}{2}, \frac{3\pi}{2}) \backslash B(0, \alpha_{j_p}, k, \mu_{j_p}), k \geq N_1\}$. Similarly, using

$$\mu_{k,j_p}(\pi + \frac{\mu_s - \mu_{j_p} \mp \alpha_{j_p}}{4n\pi k}) = (2\pi k + \pi)^n + (2\pi k)^{n-2}\frac{\mu_s + \mu_{j_p} \mp \alpha_{j_p}}{2} + O(k^{n-3}),$$

which can be proved by direct calculations, we obtain that the set

$$(h, \infty) \backslash \bigcup_{k: k \geq N_1;\ s=1,2\ldots,m} C(\pi, k, j_p, s, \alpha_{j_p}), \text{ where}$$

$C(\pi, k, j_p, s, \alpha_{j_p}) = \{x \in \mathbb{R} : |x - (2\pi k + \pi)^n + (2\pi k)^{n-2}\frac{\mu_s + \mu_{j_p}}{2}| < \alpha_{j_p}(2\pi k)^{n-2}\}$, is a subset of

$$\{\mu_{k,j_p}(t) : t \in [-\frac{\pi}{2}, \frac{3\pi}{2}) \backslash B(\pi, \alpha_{j_p}, k, \mu_{j_p}), k \geq N_1\}.$$

Now using (57) and the continuity of $\lambda_{k,j_p}(t)$ on $[-\frac{\pi}{2}, \frac{3\pi}{2}) \backslash B(\alpha_{j_p}, k, \mu_{j_p})$, we see that the set

$$(H, \infty) \backslash ((\bigcup_{k: k \geq N_1;\ s=1,2\ldots,m} C(k, j_p, s, 2\alpha_{j_p})),$$

where $H = h+1$, $C(k, j_p, s, 2\alpha_{j_p}) = C(0, k, j_p, s, 2\alpha_{j_p}) \cup C(\pi, k, j_p, s, 2\alpha_{j_p})$, is a subset of the set $\{\lambda_{k,j_p}(t) : t \in [-\frac{\pi}{2}, \frac{3\pi}{2}) \backslash B(\alpha_{j_p}, k, \mu_{j_p}), k \geq N_1\}$. Thus we have

$$\bigcup_{p=1,2,3} ((H, \infty) \backslash (\bigcup_{k \geq N_1; s=1,2,\ldots,m} C(k, j_p, s, 2\alpha_{j_p}))) \subset \sigma(L(P)).$$

To prove the inclusion $(H, \infty) \subset \sigma(L(P))$ it is enough to show that the set

$$\bigcap_{p=1,2,3} (\bigcup_{k \geq N_1; s=1,2,\ldots,m} C(k, j_p, s, 2\alpha_{j_p}))$$

is empty. If this set contains an element $x$, then

$$x \in \bigcup_{k \geq N_1; s=1,2,\ldots,m} C(k, j_p, s, \alpha_{j_p}) \qquad (60)$$

for all $p = 1, 2, 3$. Using this and the definition of $C(k, j_p, s, \alpha_{j_p})$, we obtain that there exist

$k \geq N_1$; $\nu = 0, 1$ and $s = i_p$ such that

$$\mid x - (\pi(2k+\nu))^n - \frac{\mu_{j_p} + \mu_{i_p}}{2}(2\pi k)^{n-2} \mid < 2\alpha_{j_p}(2\pi k)^{n-2}$$

for all $p = 1, 2, 3$ and hence

$$\mid \frac{\mu_{j_q} + \mu_{i_q}}{2} - \frac{\mu_{j_p} + \mu_{i_p}}{2} \mid < 4\alpha_{j_p} \tag{61}$$

for all $p, q = 1, 2, 3$. Clearly, the constant $\alpha_{j_p}$ can be chosen so that

$$8\alpha_{j_p} < \min_{i_1, i_2, i_3}(diam(\{\mu_{j_1} + \mu_{i_1},\ \mu_{j_2} + \mu_{i_2},\ \mu_{j_3} + \mu_{i_3}\})), \tag{62}$$

since, by assumption of the theorem, the right-hand side of (62) is a positive constant. If (62) holds then (61) and hence (60) does not holds which implies that $(H, \infty) \subset \sigma(L(P))$. Hence the number of the gaps in the spectrum of $L(P)$ is finite ∎